# Existence of Weak Solutions for Non-Simple Elastic Surface Models

Timothy J. Healey[1]

*Dedicated to Millard Beatty on the occasion of his 90th Birthday*

**Abstract** We consider a class of models for nonlinearly elastic surfaces in this work. We have in mind thin, highly deformable structures modeled directly as two-dimensional nonlinearly elastic continua, accounting for finite membrane and bending strains and thickness change. We assume that the stored-energy density is polyconvex with respect to the second gradient of the deformation, and we require that it grow unboundedly as the local area ratio approaches zero. For sufficiently fast growth, we show that the latter is uniformly bounded away from zero at an energy minimizer. With this in hand, we rigorously derive the weak form of the Euler-Lagrange equilibrium equations.

**Key Words** Nonlinear elasticity, thin structures, second gradient, polyconvexity, Euler-Lagrange equation

**Mathematics Subject Classification** 74B20, 35D05, 49K20

## 1. Introduction

We consider a class of models for nonlinearly elastic surfaces in this work. We have in mind thin, highly deformable structures modeled *directly* as two-dimensional nonlinearly elastic continua, accounting for finite membrane and bending strains and thickness change. As discussed by Steigmann and Ogden [17], such systems were first postulated in [4], [6]. Hilgers and Pipkin [11], [12] provide formal asymptotic derivations of such models via small-thickness expansions from bulk nonlinear elasticity. Nonlinearly elastic generalizations of such to intrinsically curved surfaces are pursued in [17]. Each of these works account for the finite-strain measures discussed above. The model treated here falls within that general class, although we assume a flat (but not-necessarily stress-free) reference configuration. Elastic surface theories that include second-gradient terms are called "non-simple" in [17]; we adopt that nomenclature here. Our main purpose is to establish the existence of weak solutions of the Euler-Lagrange equilibrium equations for a general class of mixed boundary-value problems.

An important motivation for the present work comes from our recent studies on wrinkling of highly stretchable elastomer membranes [9], [13]. The incorporation of finite nonlinear elasticity in the membrane portion of the model is crucial for capturing the correct phenomena for that class of problems. The nonlinear membrane model employed in [13] is obtained by viewing the thin structure as an incompressible, Mooney-Rivlin solid in the absence of through-thickness strain variation. On eliminating the pressure (assuming stress-free upper and lower faces), one obtains the following equivalent two-dimensional membrane energy [14]:

$$W_m = \alpha\left(tr\mathbf{C} + (\det \mathbf{C})^{-1} - 3\right) + \beta\left(\det \mathbf{C} + tr\mathbf{C}(\det \mathbf{C})^{-1} - 3\right), \tag{1.1}$$

where $\mathbf{C}$ denotes right Cauchy-Green strain 2-tensor, and $\alpha, \beta > 0$ are material constants. Note that (1.1) implies

$$W_m \nearrow \infty \text{ as } J \searrow 0, \tag{1.2}$$

---

[1] Department of Mathematics, Cornell University, Ithaca, NY 14850, USA
tjh10@cornell.edu

where $J := (\det \mathbf{C})^{1/2}$ is the local area ratio. Condition (1.2) is also an important ingredient of two-dimensional (planar) nonlinear elasticity, and in Section 2 we discuss its relationship to orientation preservation in the context of three-dimensional nonlinear elasticity, cf. Remark 2.2. Accordingly, we adopt it here as a basic requirement. It is worth mentioning that (1.2) is specified (although not employed) in [15]. The bending model used in [13] is rather crude, incorporating the linearized curvature, while thickness change is not accounted for. We consider here the more general class of nonlinearly elastic models described above in the presence of general dead loading and boundary conditions.

The outline of the work is as follows. In Section 2 we present the formulation along with our mathematical hypotheses. We assume hyper-elasticity with the stored-energy density as a function of both the gradient and second gradient of the deformation. We discuss the additive decomposition of the latter into two parts – one characterizing bending and the other typically associated with change in thickness. We assume both coercivity and polyconvexity in the second-gradient argument. Our main goal is to obtain weak solutions, the construction of which requires, among other things, sufficient smoothness. In particular, we require $p > 2$ in the coercivity condition. For simplicity, we make hypotheses consistent with that goal from the outset. We finish the section presenting the general class of boundary-value problem under consideration, which includes a wide variety of mixed boundary conditions in the presence of dead loadings.

In Section 3 we present the steps leading to the existence of weak solutions. A global energy minimizer can be obtained via appropriate specialization of the general results of Ball, Currie and Olver [5] for $k$-gradient, polyconvex systems, without assuming $p > 2$. We give a proof here for completeness, working out all details in our simpler, specific setting. We note that the resulting smoothness alone (due to embedding) is not enough to ensure a rigorous first-variation at a minimizer. We further require sufficiently fast growth as in (1.2) to show that the local-area-ratio field of a minimizer is uniformly bounded away from zero. This follows from an argument presented in [10], designed for problems of second-gradient bulk nonlinear elasticity. With this in hand, the first-variation condition at a minimizer is rigorously obtained, delivering the weak form of the Euler-Lagrange equilibrium equations. We make some final remarks in Section 4.

## 2. Problem Formulation

Let $\Omega \subset \mathbb{R}^2$ denote a bounded domain with a strongly locally Lipschitz boundary $\partial \Omega$, cf. [1]. We henceforth make the identification $\mathbb{R}^2 \cong \mathbb{R}^2 \times \{0\} \subset \mathbb{R}^3$. (As is common practice, we let $\mathbb{R}^n$ denote both Euclidean point space and its translate or tangent space.) We associate $\bar{\Omega}$ with a reference configuration for a material surface in a "flat" state embedded in $\mathbb{R}^3$. We denote deformations via $\mathbf{f} : \Omega \to \mathbb{R}^3$, and the deformation gradient or total derivative of $\mathbf{f}$ at $\mathbf{x} \in \Omega$ is denoted $\mathbf{F}(\mathbf{x}) = \nabla \mathbf{f}(\mathbf{x}) \in L(\mathbb{R}^2, \mathbb{R}^3)$. Also, $J := \sqrt{\det \mathbf{C}}$ denotes the local area ratio, where $\mathbf{C} := \mathbf{F}^T \mathbf{F} \in L(\mathbb{R}^2)$ is the right Cauchy-Green strain tensor. We define

$$L^+(\mathbb{R}^2, \mathbb{R}^3) := \{\mathbf{F} \in L(\mathbb{R}^2, \mathbb{R}^3) : J > 0\}. \tag{2.1}$$

We note that the symmetric second-order tensor $\mathbf{C}$ is positive-definite if and only if $\mathbf{F} \in L^+(\mathbb{R}^2, \mathbb{R}^3)$. The second gradient of $\mathbf{f}$ at $\mathbf{x} \in \Omega$ is denoted $\mathcal{G}(\mathbf{x}) = \nabla \mathbf{F}(\mathbf{x}) = \nabla^2 \mathbf{f}(\mathbf{x})$. Relative to the standard orthonormal basis, denoted $\{\mathbf{e}_1, \mathbf{e}_2, \mathbf{e}_3\}$, we have

$$\nabla \mathbf{f} = \mathbf{f}_{,\alpha} \otimes \mathbf{e}_\alpha = (\partial f_i / \partial x_\alpha)\mathbf{e}_i \otimes \mathbf{e}_\alpha,$$
$$\nabla^2 \mathbf{f} = \mathbf{f}_{,\alpha\beta} \otimes \mathbf{e}_\alpha \otimes \mathbf{e}_\beta = (\partial^2 f_i / \partial x_\alpha \partial x_\beta)\mathbf{e}_i \otimes \mathbf{e}_\alpha \otimes \mathbf{e}_\beta, \tag{2.2}$$

with Latin indices summing from 1-3 and Greek indices from 1-2. The triple tensor product in (2.2)$_2$ is defined by

$$(\mathbf{a} \otimes \mathbf{b} \otimes \mathbf{c})\mathbf{u} := \mathbf{a} \otimes \mathbf{b}(\mathbf{c} \cdot \mathbf{u}),$$
$$(\mathbf{a} \otimes \mathbf{b} \otimes \mathbf{c})(\mathbf{u} \otimes \mathbf{v}) := \mathbf{a}(\mathbf{b} \cdot \mathbf{u})(\mathbf{c} \cdot \mathbf{v}), \tag{2.3}$$

for all $\mathbf{a} \in \mathbb{R}^3, \mathbf{b}, \mathbf{c}, \mathbf{u}, \mathbf{v} \in \mathbb{R}^2$. From (2.2)$_2$ and (2.3)$_2$, we see that the third-order tensor $\nabla^2 \mathbf{f}(\mathbf{x})$ ($\mathbf{x}$ fixed) belongs to the space $L^s\left(L(\mathbb{R}^2), \mathbb{R}^3\right)$, the latter defined as the set of all linear transformations $\mathcal{G}: L(\mathbb{R}^2) \to \mathbb{R}^3$ satisfying $\mathcal{G}\mathbf{A} = \mathcal{G}\mathbf{A}^T$ for all $\mathbf{A} \in L(\mathbb{R}^2)$. We also define

$$\mathbf{d} \cdot (\mathbf{a} \otimes \mathbf{b} \otimes \mathbf{c}) := (\mathbf{d} \cdot \mathbf{a})\mathbf{b} \otimes \mathbf{c}, \tag{2.4}$$

for all $\mathbf{a}, \mathbf{d} \in \mathbb{R}^3, \mathbf{b}, \mathbf{c} \in \mathbb{R}^2$.

**Remark 2.1** From (2.2)$_2$ and (2.3)$_2$, it follows that every $\mathcal{G} \in L^s\left(L(\mathbb{R}^2), \mathbb{R}^3\right)$ can be associated with a symmetric bilinear transformation, say, $\tilde{\mathcal{G}}: \mathbb{R}^2 \times \mathbb{R}^2 \to \mathbb{R}^3$, via $\tilde{\mathcal{G}}[\mathbf{a}, \mathbf{b}] := \mathcal{G}(\mathbf{a} \otimes \mathbf{b}) = \mathcal{G}(\mathbf{b} \otimes \mathbf{a})$ for all $\mathbf{a}, \mathbf{b} \in \mathbb{R}^2$.

We assume that the material surface possesses a stored-energy function of the form $\Psi(\mathcal{G}, \mathbf{F})$; $\Psi: L^s\left(L(\mathbb{R}^2), \mathbb{R}^3\right) \times L^+(\mathbb{R}^2, \mathbb{R}^3) \to (0, \infty)$. We refer to such a structure as a *non-simple elastic surface*, cf. [17]. The second gradient $\mathcal{G} = \nabla^2 \mathbf{f}$ captures bending as well as membrane effects. To see this, recall that $\mathbf{n} := (\mathbf{F}\mathbf{e}_1 \times \mathbf{F}\mathbf{e}_2)/J$ defines a unit-normal field on the deformed surface $\Sigma := \mathbf{f}(\Omega)$, while $\mathbf{a}_\alpha := \mathbf{f}_{,\alpha} = \mathbf{F}\mathbf{e}_\alpha, \alpha = 1, 2$, are tangent vectors on the surface. We may then write $\mathbf{F} = \nabla \mathbf{f} = \mathbf{a}_\alpha \otimes \mathbf{e}_\alpha$ and $\nabla^2 \mathbf{f} = \mathbf{a}_{\alpha,\beta} \otimes \mathbf{e}_\alpha \otimes \mathbf{e}_\beta$, cf. (2.2). From (2.2)$_2$ and (2.4), we have $\mathbf{n} \cdot \nabla^2 \mathbf{f} = (\mathbf{n} \cdot \mathbf{a}_{\alpha,\beta})\mathbf{e}_\alpha \otimes \mathbf{e}_\beta$, where $\mathbf{n} \cdot \mathbf{a}_{\alpha,\beta} = -\mathbf{n}_{,\alpha} \cdot \mathbf{a}_\beta$ are the components of the second fundamental form for $\Sigma$. Hence,

$$\mathbf{K} := \mathbf{n} \cdot \nabla^2 \mathbf{f} = (\mathbf{n} \cdot \mathbf{f}_{,\alpha\beta})\mathbf{e}_\alpha \otimes \mathbf{e}_\beta, \tag{2.5}$$

which is called the *relative curvature tensor* in [17], has the same components as the Weingarten map $\mathbf{L} = (\mathbf{n} \cdot \mathbf{f}_{,\alpha\beta})\mathbf{a}^\alpha \otimes \mathbf{a}^\beta$, where $\{\mathbf{a}^\alpha, \mathbf{a}^\beta\}$ denotes the dual tangential basis field. From (2.2)$_2$, and (2.4), we also have

$$\mathbf{a}^\gamma \cdot \nabla^2 \mathbf{f} = (\mathbf{a}^\gamma \cdot \mathbf{f}_{,\alpha\beta})\mathbf{e}_\alpha \otimes \mathbf{e}_\beta = \Gamma^\gamma_{\alpha\beta}\mathbf{e}_\alpha \otimes \mathbf{e}_\beta, \quad \gamma = 1, 2, \tag{2.6}$$

where $\mathbf{a}^\gamma \cdot \mathbf{a}_{\alpha,\beta} = \Gamma^\gamma_{\alpha\beta}$ are the usual Christoffel symbols. Finally, (2.5) and (2.6) yield the decomposition

$$\nabla^2 \mathbf{f} = K_{\alpha\beta}\mathbf{n} \otimes \mathbf{e}_\alpha \otimes \mathbf{e}_\beta + \Gamma^\gamma_{\alpha\beta}\mathbf{a}_\gamma \otimes \mathbf{e}_\alpha \otimes \mathbf{e}_\beta. \tag{2.7}$$

Let $\mathcal{T}_y$ denote tangent space at $\mathbf{y} = \mathbf{f}(\mathbf{x})$, with $\mathbf{x} \in \Omega$ fixed. Observe that $K_{\alpha\beta}\mathbf{n} \otimes \mathbf{e}_\alpha \otimes \mathbf{e}_\beta : L(\mathbb{R}^2) \to \mathcal{T}_y^\perp$, whereas $\Gamma^\gamma_{\alpha\beta}\mathbf{a}_\gamma \otimes \mathbf{e}_\alpha \otimes \mathbf{e}_\beta : L(\mathbb{R}^2) \to \mathcal{T}_y$. Clearly, the former accounts for bending, while the latter accounts for changes in membrane thickness, e.g., [11], [12]. Note that in a planar state $\mathbf{f}(\Omega) \subset \mathbb{R}^2$, we have $\mathcal{T}_y \equiv \mathbb{R}^2$ and $\mathbf{n} \equiv \mathbf{e}_3 \Rightarrow \mathbf{K} \equiv \mathbf{0}$, cf. (2.5). Hence, in this special case we have $\nabla^2 \mathbf{f} = \Gamma^\gamma_{\alpha\beta}\mathbf{a}_\gamma \otimes \mathbf{e}_\alpha \otimes \mathbf{e}_\beta = \mathbf{f}_{,\alpha\beta} \otimes \mathbf{e}_\alpha \otimes \mathbf{e}_\beta = (\partial^2 f_\gamma / \partial x_\alpha \partial x_\beta)\mathbf{e}_\gamma \otimes \mathbf{e}_\alpha \otimes \mathbf{e}_\beta$, which is the second gradient in the plane.

**Remark 2.2** Observe that $J = (\mathbf{a}_1 \times \mathbf{a}_2) \cdot \mathbf{n} > 0$ insures that $\{\mathbf{a}_1, \mathbf{a}_2, \mathbf{n}\}$ and $\{\mathbf{e}_1, \mathbf{e}_2, \mathbf{e}_3\}$ have the same orientation. Thus, (1.2) is the remnant of the usual constitutive assumption from three-dimensional nonlinear elasticity for maintaining orientation.

We assume that the stored-energy function is objective, viz.,

$$\Psi(\mathbf{Q}\mathcal{G}, \mathbf{Q}\mathbf{F}) = \Psi(\mathcal{G}, \mathbf{F}) \text{ for all } \mathbf{Q} \in SO(3): \tag{2.8}$$

Expressing $\mathcal{G} = \mathcal{G}_{i\alpha\beta}\mathbf{e}_i \otimes \mathbf{e}_\alpha \otimes \mathbf{e}_\beta$ and $\mathbf{Q} = Q_{ij}\mathbf{e}_i \otimes \mathbf{e}_j$, then $\mathbf{Q}\mathcal{G} = Q_{ij}\mathcal{G}_{j\alpha\beta}\mathbf{e}_i \otimes \mathbf{e}_\alpha \otimes \mathbf{e}_\beta \in L^s(L(\mathbb{R}^2), \mathbb{R}^3)$. Of course, (2.8) places restrictions on the dependence of $\Psi$ on its arguments, but this plays no role in what follows. Nonetheless, we mention that (2.8) is not vacuous. For example, it holds if $\Psi$ is another real-valued function of $\mathbf{C}$ and its gradient $\nabla\mathbf{C}$.

We further assume throughout that the stored-energy function satisfies the following properties:

**H1)** For $p > 2$ and $q \geq 2p/(p-2)$, there is a constant $C > 0$ such that

$$\Psi(\mathcal{G}, \mathbf{F}) \geq C[|\mathcal{G}|^p + J^{-q}] \text{ for all } (\mathcal{G}, \mathbf{F}) \in L^s\left(L(\mathbb{R}^2), \mathbb{R}^3\right) \times L^+(\mathbb{R}^2, \mathbb{R}^3), \tag{2.9}$$

where $|\mathcal{G}|^2 = \mathcal{G} \cdot \mathcal{G} := \mathcal{G}_{i\alpha\beta}\mathcal{G}_{i\alpha\beta}$.

**H2)** $\mathcal{G} \mapsto \Psi$ is polyconvex:

Let $\mathcal{D}^{[2]}$ denote the list of all $2 \times 2$ sub-determinants of the components of $\mathcal{G}$. Arranging the latter, viz., $\mathcal{G}_{i\alpha\beta}$, into a $6 \times 2$ matrix, we see that there are $\binom{6}{2} = 15$ such independent determinants, written $\mathcal{D}^{[2]} = (d_1, d_2, \ldots, d_{15})$. These are listed below, expressed in terms of $\mathcal{G}_{i\alpha\beta} = f_{i,\alpha\beta}$:

$$\begin{aligned} d_\ell[\mathbf{f}] &= f_{\ell,11}f_{\ell,22} - (f_{\ell,12})^2, \quad \ell = 1,2,3, \\ &= f_{1,11}f_{(\ell-2),12} - f_{1,12}f_{(\ell-2),11}, \quad \ell = 4,5, \\ &= f_{1,11}f_{(\ell-4),22} - f_{1,12}f_{(\ell-4),12}, \quad \ell = 6,7, \\ &= f_{1,12}f_{(\ell-6),12} - f_{1,22}f_{(\ell-6),11}, \quad \ell = 8,9, \\ &= f_{1,12}f_{(\ell-8),22} - f_{1,22}f_{(\ell-8),12}, \quad \ell = 10,11, \\ &= f_{2,11}f_{3,(\ell-11)2} - f_{2,12}f_{3,1(\ell-11)}, \quad \ell = 12,13, \\ &= f_{2,12}f_{3,(\ell-13)2} - f_{2,22}f_{3,1(\ell-13)}, \quad \ell = 14,15. \end{aligned} \tag{2.10}$$

A function $\varphi: L^s(L(\mathbb{R}^2), \mathbb{R}^3) \to \mathbb{R}$ is said to be *polyconvex* if there is a convex function $\tilde{\varphi}: L^s(L(\mathbb{R}^2), \mathbb{R}^3) \times \mathbb{R}^{15} \to \mathbb{R}$ such that $\varphi(\mathcal{G}) = \tilde{\varphi}(\mathcal{G}, \mathcal{D}^{[2]})$, cf. [5]. Here we write $\Psi(\mathcal{G}, \mathbf{F}) \equiv \Phi(\mathcal{G}, \mathcal{D}^{[2]}, \mathbf{F})$, with $(\mathcal{G}, \mathcal{D}^{[2]}) \mapsto \Phi$ convex, and we assume that $\Phi: L^s(L(\mathbb{R}^2), \mathbb{R}^3) \times \mathbb{R}^{15} \times L^+(\mathbb{R}^2, \mathbb{R}^3) \to (0, \infty)$ is $C^1$.

**H3)** For $p > 2$, $J \geq \eta > 0$, and $|\mathbf{F}|^2 = \mathbf{F} \cdot \mathbf{F} := tr\mathbf{C} \leq R^2$, there is a constant $C_{\eta, R} > 0$ such that

$$\Psi(\mathcal{G}, \mathbf{F}) \leq C_{\eta, R}(|\mathcal{G}|^p + 1),$$
$$|\Psi_{\mathcal{G}}(\mathcal{G}, \mathbf{F})| \leq C_{\eta, R}(|\mathcal{G}|^{p-1} + 1), \quad |\Psi_{\mathbf{F}}(\mathcal{G}, \mathbf{F})| \leq C_{\eta, R}(|\mathcal{G}|^p + 1). \tag{2.11}$$

We consider a class of mixed boundary-value problems as follows. We let $W^{k,p}(\Omega, \mathbb{R}^3)$ denote the Sobolev space of vector-valued $p$-integrable functions, such that all weak partial derivatives of order less than or equal to $k$ are also $p$-integrable. Here we are interested in $p > 2$ and $k = 0, 1, 2$, with $W^{0,p}(\Omega, \mathbb{R}^3) \equiv L^p(\Omega, \mathbb{R}^3)$. The norms are defined by

$$\|\mathbf{f}\|^p_{L^p(\Omega, \mathbb{R}^3)} = \int_\Omega |\mathbf{f}|^p \, dx,$$
$$\|\mathbf{f}\|^p_{W^{1,p}(\Omega, \mathbb{R}^3)} = \|\mathbf{f}\|^p_{L^p(\Omega, \mathbb{R}^3)} + \int_\Omega |\nabla \mathbf{f}|^p \, dx, \tag{2.12}$$
$$\|\mathbf{f}\|^p_{W^{2,p}(\Omega, \mathbb{R}^3)} = \|\mathbf{f}\|^p_{W^{1,p}(\Omega, \mathbb{R}^3)} + \int_\Omega |\nabla^2 \mathbf{f}|^p \, dx,$$

respectively, where the tensor-Euclidean norms are employed in the integrands in (2.12)$_{2,3}$, as defined in (H1), (H3). Consider a subset $\Gamma \subset \partial\Omega$ with positive length, and let

$$W^{2,p}_\Gamma(\Omega, \mathbb{R}^3) = \{\mathbf{u} \in W^{2,p}(\Omega, \mathbb{R}^3) : \mathbf{u} = \mathbf{0} \text{ and } [\nabla \mathbf{u}]\mathbf{v} = \mathbf{0} \text{ a.e. on } \Gamma\}, \tag{2.13}$$

where $\mathbf{v} \in \mathbb{R}^2$ denotes the outward unit normal field, which exists *a.e.* on $\partial\Omega$. The usual Poincaré inequality gives $\|\mathbf{u}\|_{L^p(\Omega, \mathbb{R}^3)} \leq C_1 \|\nabla \mathbf{u}\|_{L^p(\Omega, \mathbb{R}^3)} \leq C_2 \|\nabla^2 \mathbf{u}\|_{L^p(\Omega, \mathbb{R}^3)}$ for all $\mathbf{u} \in W^{2,p}_\Gamma(\Omega, \mathbb{R}^3)$. Thus, we have

$$\|\mathbf{u}\|_{W^{2,p}(\Omega, \mathbb{R}^3)} \leq C \|\nabla^2 \mathbf{u}\|_{L^p(\Omega, \mathbb{R}^3)}, \tag{2.14}$$

for all $\mathbf{u} \in W^{2,p}_\Gamma(\Omega, \mathbb{R}^3)$, i.e., the right side of (2.14) defines an equivalent norm.

We consider the admissible set

$$\mathcal{A} = \{\mathbf{f} \in W^{2,p}(\Omega, \mathbb{R}^3) : \mathbf{f} - \mathbf{f}_o \in W^{2,p}_\Gamma(\Omega, \mathbb{R}^3), \nabla \mathbf{f} \in L^+(\mathbb{R}^2, \mathbb{R}^3) \text{ a.e. in } \Omega\}, \tag{2.15}$$

where $\mathbf{f}_o \in W^{2,p}(\Omega, \mathbb{R}^3)$ is prescribed. The following potential-energy functional, $E: \mathcal{A} \to \mathbb{R}$, is to be minimized:

$$E[\mathbf{f}] := \int_\Omega \left[ \Psi(\nabla^2 \mathbf{f}, \nabla \mathbf{f}) - (\mathbf{b} \cdot \mathbf{f} + \mathbf{B} \cdot \nabla \mathbf{f}) \right] dx - \int_{\Gamma_c} \left[ \boldsymbol{\tau} \cdot \mathbf{f} + \boldsymbol{\mu} \cdot ([\nabla \mathbf{f}]\mathbf{v}) \right] ds, \tag{2.16}$$

where $\mathbf{b} \in L^1(\Omega, \mathbb{R}^3), \mathbf{B} \in L^1(\Omega, L(\mathbb{R}^2, \mathbb{R}^3))$ represent prescribed body-force and generalized body-force densities, respectively, $\mathbf{B} \cdot \mathbf{F} := B_{i\alpha} F_{i\alpha}$, $\boldsymbol{\tau}, \boldsymbol{\mu} \in L^1(\Gamma_c, \mathbb{R}^3)$ are prescribed surface-traction and surface-hyper-traction densities, respectively, and $\Gamma_c := \partial\Omega \setminus \Gamma$.

We observe that in the special case $\Gamma \equiv \partial\Omega$, i.e., "clamped" Dirichlet conditions are specified in (2.13), $W^{2,p}_o(\Omega, \mathbb{R}^3)$ replaces $W^{2,p}_\Gamma(\Omega, \mathbb{R}^3)$ in $\mathcal{A}$, and the surface integral is not present in (2.16). It's worth noting that no smoothness conditions on $\partial\Omega$ are required in this case. Weaker Dirichlet conditions, viz., "pinned" conditions are also of interest. The potential energy functional to be minimized in this case is given by

$$\tilde{E}[\mathbf{f}] := \int_\Omega \left[\Psi(\nabla^2\mathbf{f}, \nabla\mathbf{f}) - (\mathbf{b}\cdot\mathbf{f} + \mathbf{B}\cdot\nabla\mathbf{f})\right]dx - \int_{\partial\Omega}\boldsymbol{\mu}\cdot([\nabla\mathbf{f}]\mathbf{v})ds, \tag{2.17}$$

over the admissible set

$$\tilde{\mathcal{A}} = \{\mathbf{f}\in W^{2,p}(\Omega,\mathbb{R}^3) : \mathbf{f} = \mathbf{f}_o \text{ on } \partial\Omega, \nabla\mathbf{f}\in L^+(\mathbb{R}^2,\mathbb{R}^3) \text{ a.e. in } \Omega\}. \tag{2.18}$$

We claim that (2.14) is also valid for all elements belonging to

$$W^{2,p}_{\partial\Omega}(\Omega,\mathbb{R}^3) = \{\mathbf{u}\in W^{2,p}(\Omega,\mathbb{R}^3) : \mathbf{u}=\mathbf{0} \text{ on } \partial\Omega\}. \tag{2.19}$$

As before, $\|\mathbf{u}\|_{L^p(\Omega,\mathbb{R}^3)} \leq C_1\|\nabla\mathbf{u}\|_{L^p(\Omega,\mathbb{R}^3)}$ for all $\mathbf{u}\in W^{2,p}_{\partial\Omega}(\Omega,\mathbb{R}^3)$, while a generalized Poincaré inequality reads

$$\|\nabla\mathbf{u}\|^p_{L^p(\Omega,\mathbb{R}^3)} \leq C_2\left[\|\nabla^2\mathbf{u}\|^p_{L^p(\Omega,\mathbb{R}^3)} + \left|\int_\Omega \nabla\mathbf{u}\,dx\right|^p\right].$$

But each component of the integral in the second term above vanishes by virtue of Green's theorem for all $\mathbf{u}\in W^{2,p}_{\partial\Omega}(\Omega,\mathbb{R}^3)$, cf. [16].

As a final remark we note that, due to embedding $(p > 2)$, the boundary-value prescriptions in (2.15) can be interpreted in the pointwise sense.

## 3. Energy Minima and Weak Solutions

We prove the existence of weak solutions in this section, culminating in Theorem 3.4. Throughout we presume hypotheses (H1) – (H3). The first step is:

**Proposition 3.1.** Assume that $\mathbf{f}_o \in W^{2,p}(\Omega,\mathbb{R}^3)$ and $\nabla\mathbf{f}_o \in L^+(\mathbb{R}^2,\mathbb{R}^3)$ on $\bar{\Omega}$. Then $E$ attains its minimum on $\mathcal{A}$, viz., there exists $\mathbf{f}_* \in \mathcal{A}$ such that $E[\mathbf{f}_*] = \inf_{\mathbf{f}\in\mathcal{A}} E[\mathbf{f}]$. Likewise, there exists $\mathbf{f}_* \in \tilde{\mathcal{A}}$ such that $\tilde{E}[\mathbf{f}_*] = \inf_{\mathbf{f}\in\tilde{\mathcal{A}}} \tilde{E}[\mathbf{f}]$.

**Proof:** By embedding, $\mathbf{f}_o \in C^1(\bar{\Omega},\mathbb{R}^3)$, and by assumption, $J_o = (\det[\nabla\mathbf{f}_o^T\nabla\mathbf{f}_o])^{1/2} \geq m$ on $\bar{\Omega}$, where $m > 0$ denotes its minimum on the compact set. Then from $(2.11)_1$ and (2.16), we see that $E[\mathbf{f}_o] < \infty$. For any $\mathbf{f}\in\mathcal{A}$, it follows from (H1) and (2.16) that

$$\begin{aligned}C\|\nabla^2\mathbf{f}\|^p_{L^p(\Omega,\mathbb{R}^3)} &\leq E[\mathbf{f}] + \max_{\bar{\Omega}}|\mathbf{f}|\left(\|\mathbf{b}\|_{L^1(\Omega,\mathbb{R}^3)} + \|\boldsymbol{\tau}\|_{L^1(\Gamma_c,\mathbb{R}^3)}\right) \\ &\quad + \max_{\bar{\Omega}}|\nabla\mathbf{f}|\left(\|\mathbf{B}\|_{L^1(\Omega,L(\mathbb{R}^3))} + \|\boldsymbol{\mu}\|_{L^1(\Gamma_c,\mathbb{R}^3)}\right) \\ &\leq E[\mathbf{f}] + M\|\mathbf{f}\|_{C^1(\Omega,\mathbb{R}^3)} \leq E[\mathbf{f}] + \tilde{M}\|\mathbf{f}\|_{W^{2,p}(\Omega,\mathbb{R}^3)},\end{aligned} \tag{3.1}$$

the last inequality of which follows by embedding. Since $\mathbf{f} - \mathbf{f}_o \in W^{2,p}_\Gamma(\Omega,\mathbb{R}^3)$, inequality (2.14) implies

$$\|\mathbf{f}\|_{W^{2,p}(\Omega,\mathbb{R}^3)} - \|\mathbf{f}_o\|_{W^{2,p}(\Omega,\mathbb{R}^3)} \leq \tilde{C}\left(\|\nabla^2\mathbf{f}\|_{L^p(\Omega,\mathbb{R}^3)} + \|\nabla^2\mathbf{f}_o\|_{L^p(\Omega,\mathbb{R}^3)}\right),$$

which combined with (3.1) yields

$$E[\mathbf{f}] \geq C_1 \|\mathbf{f}\|^p_{W^{2,p}(\Omega,\mathbb{R}^3)} + C_2, \qquad (3.2)$$

for constants $C_1 > 0$ and $C_2$. Let $\{\mathbf{f}_j\} \subset \mathcal{A}$ be an infimizing sequence. Since $\inf_{\mathbf{f} \in \mathcal{A}} E[\mathbf{f}] \leq E[\mathbf{f}_o] < \infty$, inequality (3.2) implies that $\{\mathbf{f}_j\} \subset W^{2,p}(\Omega)$ is uniformly bounded. Hence, there exists a weakly convergent subsequence,

$$\mathbf{f}_{j_k} \rightharpoonup \mathbf{f}_* \text{ weakly in } W^{2,p}(\Omega,\mathbb{R}^3), \text{ and } \mathbf{f}_{j_k} \to \mathbf{f}_* \text{ strongly in } C^1(\bar{\Omega},\mathbb{R}^3), \qquad (3.3)$$

the latter of which follows by compact embedding. The closed linear subspace $W^{2,p}_\Gamma(\Omega,\mathbb{R}^3) \subset W^{2,p}(\Omega,\mathbb{R}^3)$ is also weakly closed, and therefore $\mathbf{f}_* - \mathbf{f}_o \in W^{2,p}_\Gamma(\Omega,\mathbb{R}^3)$. Also, from (2.15) and (3.3) we deduce that $J_{j_k} := (\det[\nabla \mathbf{f}^T_{j_k} \nabla \mathbf{f}_{j_k}])^{1/2}$ converges uniformly to $J_* \geq 0$ on $\bar{\Omega}$. We claim that $J_* > 0$ a.e. in $\Omega$. If not, then $J_{j_k} \to 0$ uniformly in $\Upsilon \subset \Omega$, where $|\Upsilon| > 0$. By virtue of (H1) and Fatou's lemma, we have $\liminf_{k \to \infty} \int_\Upsilon \Psi(\nabla^2 \mathbf{f}_{j_k}(\mathbf{x}), \nabla \mathbf{f}_{j_k}(\mathbf{x})) dx \geq \int_\Upsilon \liminf_{k \to \infty} \Psi(\nabla^2 \mathbf{f}_{j_k}(\mathbf{x}), \nabla \mathbf{f}_{j_k}(\mathbf{x})) dx \to \infty$. But in view of (2.16), this contradicts the fact that $\lim_{k \to \infty} E[\mathbf{f}_{j_k}] = \inf_{\mathbf{f} \in \mathcal{A}} E[\mathbf{f}] < \infty$. We conclude that $\mathbf{f}_* \in \mathcal{A}$.

The weak convergence of $\{\mathbf{f}_{j_k}\} \subset W^{2,p}(\Omega)$ implies that each of the sub-determinants in (2.10) converges weakly, viz.,

$$d_\ell[\mathbf{f}_{j_k}] \rightharpoonup d_\ell[\mathbf{f}_*] \text{ weakly in } L^{p/2}(\Omega), \text{ for } \ell = 1, 2, \ldots, 15. \qquad (3.4)$$

For example, define the vector field $\mathbf{w}(\mathbf{x}) := (f_{2,1}(\mathbf{x}), f_{3,1}(\mathbf{x}))$. Then the total derivative is given by $D\mathbf{w} = \begin{bmatrix} f_{2,11} & f_{2,12} \\ f_{3,11} & f_{3,12} \end{bmatrix}$, and $d_{12}[\mathbf{f}] = \det D\mathbf{w}$. Now consider the sequence $\{d_{12}[\mathbf{f}_{j_k}] = \det D\mathbf{w}_{j_k}\}$, which is bounded in $L^{p/2}(\Omega)$. That the weak limit is $d_{12}[\mathbf{f}_*]$ as indicated in (3.4), follows from the fact that the determinant can be written as a divergence, ultimately leading to

$$\int_\Omega \det[D\mathbf{w}_{j_k}] \varphi \, dx = -\frac{1}{2} \int_\Omega \left( Cof[D\mathbf{w}_{j_k}] \nabla \varphi \right) \cdot \mathbf{w}_{j_k} dx, \qquad (3.5)$$

for all smooth test functions $\varphi$ with compact support in $\Omega$, cf. [7], [8]. In our setting here, we have $Cof[D\mathbf{w}] = \begin{bmatrix} f_{3,12} & -f_{3,11} \\ -f_{2,12} & f_{2,11} \end{bmatrix}$. The convergence properties (3.3) then enable taking a rigorous limit as $k \to \infty$ on the right side of (3.5), which yields (3.4).

With this in hand, polyconvexity (H2) leads to the weak lower semi-continuity of $E[\cdot]$, i.e., $\liminf_{k \to \infty} E[\mathbf{f}_k] \geq E[\mathbf{f}]$ whenever $\mathbf{f}_k \rightharpoonup \mathbf{f}$ weakly in $W^{2,p}(\Omega,\mathbb{R}^3)$. First, as indicated by the estimate (3.1), the loading terms in (2.16) define a bounded linear functional, viz.,

$$\ell[\mathbf{f}] := \int_\Omega (\mathbf{b} \cdot \mathbf{f} + \mathbf{B} \cdot \nabla \mathbf{f}) dx + \int_{\Gamma_c} [\boldsymbol{\tau} \cdot \mathbf{f} + \boldsymbol{\mu} \cdot ([\nabla \mathbf{f}] \mathbf{v})] ds, \qquad (3.6)$$

with $|\ell[\mathbf{f}]| \leq M \|\mathbf{f}\|_{C^1(\bar{\Omega},\mathbb{R}^3)}$ for all $\mathbf{f} \in W^{2,p}(\Omega,\mathbb{R}^3)$, and in view of (3.3), continuity is clear. Hence, only the internal potential energy requires attention. From (H2), we deduce

$$\int_\Omega \Psi(\nabla^2 \mathbf{f}_k, \nabla \mathbf{f}_k)dx \geq \int_{\Omega_\varepsilon} \Psi(\nabla^2 \mathbf{f}_k, \nabla \mathbf{f}_k)dx$$

$$\geq \int_{\Omega_\varepsilon} \Phi(\nabla^2 \mathbf{f}, d_1[\mathbf{f}],\ldots,d_{15}[\mathbf{f}], \nabla \mathbf{f}_k)dx$$

$$+ \int_{\Omega_\varepsilon} D_1 \Phi(\nabla^2 \mathbf{f}, d_1[\mathbf{f}],\ldots,d_{15}[\mathbf{f}], \nabla \mathbf{f}_k) \cdot (\nabla^2 \mathbf{f}_k - \nabla^2 \mathbf{f})dx \quad (3.7)$$

$$+ \int_{\Omega_\varepsilon} \sum_{\ell=1}^{15} D_{d_\ell} \Phi(\nabla^2 \mathbf{f}, d_1[\mathbf{f}],\ldots,d_{15}[\mathbf{f}], \nabla \mathbf{f}_k)(d_\ell[\mathbf{f}_k] - d_\ell[\mathbf{f}])dx,$$

where $\Omega_\varepsilon := \{\mathbf{x} \in \Omega : |\nabla^2 \mathbf{f}(\mathbf{x})| \leq 1/\varepsilon,\ J(\mathbf{x}) \geq \varepsilon > 0\}$. By virtue of (3.3) and (3.4), each of the last two integrals on the right side of (3.7) approaches zero, and the first converges to $\int_{\Omega_\varepsilon} \Phi(\nabla^2 \mathbf{f}, d_1[\mathbf{f}],\ldots,d_{15}[\mathbf{f}], \nabla \mathbf{f})dA$ in the limit as $k \to \infty$. Since this holds true for all $\varepsilon > 0$, with $|\Omega - \Omega_\varepsilon| \to 0$ as $\varepsilon \to 0$, weak lower semi-continuity then follows from the Lebesgue monotone convergence theorem, cf. [7], [8]. Finally, since $E[\mathbf{f}_*] \leq \liminf_{k \to \infty} E[\mathbf{f}_{j_k}]$ and $\mathbf{f}_* \in \mathcal{A}$, the proof for $E[\cdot]$ on $\mathcal{A}$ is complete. The proof that $\tilde{E}[\mathbf{f}]$ attains its infimum on $\tilde{\mathcal{A}}$ is the same. □

**Remark 3.2.** The assumptions on $\mathbf{f}_o$, insuring that $\inf_{\mathbf{f} \in \mathcal{A}} E[\mathbf{f}] < \infty$, are made here for convenience; otherwise, we may simply assume that $\inf_{\mathbf{f} \in \mathcal{A}} E[\mathbf{f}] < \infty$.

As in the case of bulk 2D nonlinear elasticity, it is conceivable that the minimum energy configuration $\mathbf{f}_* \in \mathcal{A}$ is characterized by $J_*^2 = \det[\nabla \mathbf{f}_*^T \nabla \mathbf{f}_*] = 0$ on some set of measure zero within $\Omega$, in which case $\mathbf{x} \mapsto \Psi(\nabla^2 \mathbf{f}_*, \nabla \mathbf{f}_*)$ is infinite on that same set, cf. (H1). Of course, this is a major impediment to taking the first variation rigorously. However, with (H1) in hand, we can use a construction from [10], designed for problems in second-gradient bulk nonlinear elasticity, to show that $J_*$ is strictly positive on $\overline{\Omega}$.

**Lemma 3.3.** Given the hypotheses of Proposition 3.1, there is a number $\eta > 0$ such that $J_* = (\det[\nabla \mathbf{f}_*^T \nabla \mathbf{f}_*])^{1/2} \geq \eta$ on $\overline{\Omega}$.

**Proof.** By embedding, a minimizer $\mathbf{f}_*$ belongs to the Hölder space $C^{1,\alpha}(\overline{\Omega}, \mathbb{R}^3)$, where $\alpha = 1 - 2/p$. It follows that $J_* \in C^\alpha(\overline{\Omega})$; there is a constant $M > 0$ such that

$$J_*(\mathbf{x}) \leq J_*(\mathbf{y}) + M|\mathbf{x} - \mathbf{y}|^\alpha \quad \text{for all } \mathbf{x}, \mathbf{y} \in \overline{\Omega}. \quad (3.8)$$

Since $\partial \Omega$ is locally Lipschitz, we know that each $\mathbf{x} \in \Omega$ is the vertex of an open cone, $V_\delta^I(\mathbf{x}) \subset \Omega$, where $\delta > 0$ denotes the radius and $I \subset (0, 2\pi)$ is the domain of the polar angle, cf. [1]. Consider the integral

$$\int_{V_\delta^I(\mathbf{y})} [J_*(\mathbf{y}) + M|x - y|^\alpha]^{-q} dx = \gamma \int_0^\delta [J_*(\mathbf{y}) + Mr^\alpha]^{-q} r\, dr. \quad (3.9)$$

where $r = |\mathbf{x} - \mathbf{y}|$, $\gamma = |I|$, and $q$ is specified in (H1). Observe that the function $h:(0,\infty) \to (0,\infty)$, defined by

$$h(t) := \gamma \int_0^\delta [t + Mr^\alpha]^{-q} r\, dr, \quad (3.10)$$

is monotonically decreasing, with $h \nearrow \infty$ as $t \searrow 0$, the latter property following from (H1), viz., $-q\alpha + 1 = -q[(p-2)/p] + 1 \leq -1$. Hence, there is a constant $\eta > 0$ such that

$$h(t) \leq L_* \Leftrightarrow t \geq \eta. \tag{3.11}$$

For any $y \in \Omega$, there is a cone $V_\delta^1(y) \subset \Omega \cap D_\delta(y)$, where $D_\delta(y)$ denotes the open disk centered at $y$. From (3.6) and (3.8)-(3.10), we then deduce

$$\begin{aligned} h(J(y)) &\leq \int_{\Omega \cap D_\delta(y)} [J_*(\mathbf{y}) + M|x-y|^\alpha]^{-q} dx \\ &\leq \int_{\Omega \cap D_\delta(y)} [J_*(\mathbf{x})]^{-q} dx \leq C\big(E[\mathbf{f}_*] + \ell[\mathbf{f}_*]\big) := L_*, \end{aligned} \tag{3.12}$$

and the result follows from (3.11), which extends to all $y \in \bar{\Omega}$ via continuity. □

With Lemma 3.3 in hand, we obtain the main result:

**Theorem 3.4.** Given the hypotheses of Proposition 3.1, the minimizer $\mathbf{f}_*$ satisfies the weak form of Euler-Lagrange equilibrium equations, viz., for problem (2.16) we obtain

$$\begin{aligned} &\int_\Omega [\Psi_{\mathcal{G}}(\nabla^2 \mathbf{f}_*, \nabla \mathbf{f}_*) \cdot \nabla^2 \boldsymbol{\varphi} + \Psi_{\mathbf{F}}(\nabla^2 \mathbf{f}_*, \nabla \mathbf{f}_*) \cdot \nabla \boldsymbol{\varphi}] dx \\ &- \int_\Omega (\mathbf{b} \cdot \boldsymbol{\varphi} + \mathbf{B} \cdot \nabla \boldsymbol{\varphi}) dx - \int_{\Gamma_c} [\boldsymbol{\tau} \cdot \boldsymbol{\varphi} + \boldsymbol{\mu} \cdot ([\nabla \boldsymbol{\varphi}]\mathbf{v})] ds = 0, \end{aligned} \tag{3.13}$$

for all $\boldsymbol{\varphi} \in W_\Gamma^{2,p}(\Omega, \mathbb{R}^3)$, cf. (2.13). Likewise, (3.13) holds with $\boldsymbol{\mu} \equiv 0$ for problem (2.17) for all $\boldsymbol{\varphi} \in W_{\partial \Omega}^{2,p}(\Omega, \mathbb{R}^3)$, cf. (2.19).

**Proof.** The Gâteaux differentiability of the loading functional (3.6) is clear; we focus on the internal energy. Defining $i(t) := \int_\Omega \Psi(\nabla^2 \mathbf{f}_* + t\nabla^2 \boldsymbol{\varphi}, \nabla \mathbf{f}_* + t\nabla \boldsymbol{\varphi}) dx$, we need to rigorously justify the computation $i'(0) = \lim_{t \to 0}\{[i(t) - i(0)]/t\}$, giving the first line of (3.13). For $t \neq 0$, the fundamental theorem of calculus yields

$$[i(t) - i(0)]/t = \int_\Omega H_t(\mathbf{x}) dx, \tag{3.14}$$

where

$$H_t(\mathbf{x}) := \int_0^1 [\Psi_{\mathcal{G}}(\nabla^2 \mathbf{f}_* + st\nabla^2 \boldsymbol{\varphi}, \nabla \mathbf{f}_* + st\nabla \boldsymbol{\varphi}) \cdot \nabla^2 \boldsymbol{\varphi} + \Psi_{\mathbf{F}}(\nabla^2 \mathbf{f}_* + st\nabla^2 \boldsymbol{\varphi}, \nabla \mathbf{f}_* + st\nabla \boldsymbol{\varphi}) \cdot \nabla \boldsymbol{\varphi}] ds.$$

By embedding, $\mathbf{f}_*, \boldsymbol{\varphi} \in C^1(\bar{\Omega}, \mathbb{R}^3)$, and thus $\max_{\substack{\mathbf{x} \in \bar{\Omega}, s \in [0,1], \\ t \in [-1,1]}} |\nabla \mathbf{f}_* + st\nabla \boldsymbol{\varphi}| \leq R$. With this and Lemma 3.3 in hand, we use (2.11)$_{2,3}$ of (H3) to deduce

$$|H_t| \leq C_{\eta, R} \int_0^1 \{[|\nabla^2 \mathbf{f}_* + st\nabla^2 \boldsymbol{\varphi}|^{p-1} + 1]|\nabla^2 \boldsymbol{\varphi}| + [|\nabla^2 \mathbf{f}_* + st\nabla^2 \boldsymbol{\varphi}|^p + 1]|\nabla \boldsymbol{\varphi}|\} ds. \tag{3.15}$$

Clearly, the second term of the above integrand belongs to $L^1(\Omega)$, and an application of Hölder's inequality shows that the same holds for the first term as well, for any $(s,t) \in [0,1] \times [-1,1]$. Let $\gamma(\mathbf{x})$ denote the supremum of the sum of these two terms in $(s,t) \in [0,1] \times [-1,1]$, which gives $|H_t(\mathbf{x})| \leq \gamma(\mathbf{x})$ a.e. in $\Omega$, with $\gamma \in L^1(\Omega)$. Hence, the desired result follows by taking the limit in (3.14) as $t \to 0$ via the Lebesgue dominated convergence theorem, cf. [7], [8]. □

## 4. Concluding Remarks.

Polyconvexity in the second-gradient argument (H2), although quite general and mathematically expedient, does not easily lend itself to direct physical interpretation. Of course, convexity in the second-gradient argument alone is a special, physically reasonable case. Also, as mentioned in [5], the Gaussian curvature of the deformed surface, $\kappa := \det \mathbf{L}$, is a linear function of some of the determinants given in (2.10). Indeed, a formula from differential geometry gives $\det \mathbf{K} = J^2 \kappa$. Then writing $K_{\alpha\beta} = n_i f_{i,\alpha\beta}$, where $n_i, f_i, i = 1, 2, 3$, denote the Cartesian components, equations (2.5) and (2.10) yield

$$\kappa = [n_i^2 d_i + n_1 n_2 (d_6 - d_8) + n_2 n_3 (d_{13} - d_{14}) + n_1 n_3 (d_7 - d_9)] / J^2.$$

Even with $p > 2$, the proof of Lemma 3.3 is no longer valid if we simply require $\Psi \nearrow \infty$ as $J \searrow 0$ without the specific blow-up condition in (2.9); the existence of a weak solution is then an open question. This is also the case for bulk second-gradient nonlinear elasticity [10]. For instance, this occurs in surface models incorporating the growth condition in (1.1), where $q = 2 < 2p/(p-2)$ for all $p > 2$. We mention that the existence of an orientation-preserving, energy minimizer for our model with $4/3 < p \leq 2$ can be established. This requires a growth condition involving each of the fifteen quantities in (2.10) along with that in (2.9). Of course, our approach to weak solutions fails in that range, due to a lack of smoothness. If we ignore the requirement (1.2) altogether, then (H3) makes sense even if $J$ vanishes, and the existence of a weak solution follows routinely. However, as in bulk nonlinear elasticity, we then risk the possibility of solutions characterized by the interpenetration of matter.

In addition to [11], [12], there is an enormous literature – far too numerous to list – of works deriving surface (plate/shell) models from three-dimensional nonlinear elasticity via small-thickness asymptotics relative to a stress-free reference configuration. The results range from the rigorous to the entirely formal. Such derivations often deliver quadratic bending energies in the absence of any condition akin to (1.2) and without accounting for thickness change. Moreover, small (or even zero) membrane strains are usually a feature of the asymptotic model. These properties arise from the limiting behavior of the model as the thickness goes to zero, whereas real-world thin structures have small but finite thicknesses. While such dimensionally reduced models are well suited for certain applications, e.g., inextensible and nearly inextensible sheets and for local post-buckling in traditional plates and shells, they are incapable of capturing finite-strain phenomena for thin structures. For example, the inadequacy of the Föpple-von Kármán model in predicting the correct wrinkling phenomena in [13] is well documented in that work.

Local-injectivity conditions like (1.2) are advocated in [3] in the context of direct nonlinearly elastic Cosserat plate-shell theories. These models include through-thickness shear strains in addition to the finite-strain measures accounted for here in this work. A more stringent version of (1.2) is advocated in [3], also accounting for the small but finite thickness, based on considerations from three-dimensional nonlinear elasticity. We also mention that rigorous energy-minimization results enforcing such a condition (with a concomitant blow-up of energy) are presented in [2] for a class of nonlinearly elastic Cosserat shells satisfying the Kirchhoff-Love hypothesis.

Our approach to constructing weak solutions is not directly applicable to elastic plate-shell models incorporating finite membrane strains without thickness change and in the presence of (1.2), e.g., such as that employed in [13]. Here we refer to models characterized by a stored-energy density that is function of the deformation gradient and the bending strain (2.5) only. For growth conditions incorporating these strain measures, the smoothness needed to apply Lemma 3.3 is unclear; this is also the case for the results in [2].

Finally, we mention that the general class of surface models considered here, accounting for thickness change as well as the other finite-strain measures, is undoubtedly superior to that employed in [13] for the prediction of wrinkling in thin, highly stretched sheets.

## Acknowledgements

This work was supported in part by the National Science Foundation through grant DMS-2006586.

## References

1. Adams, R.A.: Sobolev Spaces. Academic Press, New York (1975)
2. Anicic, S.: Polyconvexity and existence theorems for nonlinearly elastic shells. J. Elasticity **132**, 161-173 (2018).
3. Antman, S.S.: Nonlinear Problems of Elasticity, 2$^{nd}$ Ed. Springer, New York (2005)
4. Balaban, M.M., Green, A.E., Naghdi, P.M.: Simple force multipoles in the theory of deformable surfaces. J. Math. Phys. **8**, 1026-1036 (1967)
5. Ball, J.M, Currie, J.C., Olver, P.J.: Null Lagrangians, weak continuity, and variational problems of arbitrary order. J. Funct. Anal. **41**, 135-174 (1981)
6. Cohen, H., DeSilva, C.N.: On a nonlinear theory of elastic shells. J. Mécanique **7**, 459-464 (1968)
7. Dacorogna, B.: Direct Methods in the Calculus of Variations, 2$^{nd}$ Ed. Springer, New York (2008)
8. Evans, L.C.: Partial Differential Equations, 2$^{nd}$ Ed. American Mathematical Society, Providence (2010)
9. Healey, T.J., Li, Q., Cheng, R.-B.: Wrinkling behavior of highly stretched rectangular elastic films via parametric global bifurcation. J. Nonlinear Sci. **23**, 777-805 (2013)
10. Healey, T.J., Krömer, S.: Injective weak solutions in second-gradient nonlinear elasticity. ESAIM: COCV **15**, 863-871 (2009)
11. Hilgers, M.G., Pipkin, A.C.: Bending energy of highly elastic membranes. Quart. Appl. Math. **50**, 389-400 (1992)
12. __________: Bending energy of highly elastic membranes II. Quart. Appl. Math. **54**, 307-316 (1996)
13. Li, Q., Healey, T.J.: Stability boundaries for wrinkling in highly stretches elastic sheets, J. Mech. Phys. Solids **97**, 260-274 (2016)
14. Müller, I, Strehlow, P.: Rubber and Rubber Balloons. Springer, Berlin (2004)
15. Pipkin, A.C.: Relaxed energy densities for large deformations of membranes. IMA J. Appl. Math. **52**, 297-308 (1994)
16. Rektorys, K.: Variational Methods in Mathematics, Science and Engineering, 2$^{nd}$ Ed. D. Reidel, Dordretcht (1980)
17. Steigmann, D.J., Ogden, R.W.: Elastic surface-substrate interactions. Proc. R. Soc. Lond. A **455**, 437-474 (1993)